\documentclass[12pt]{amsart}
\theoremstyle{plain}
\newtheorem{Thm}{Theorem}

\newtheorem{Prop}[Thm]{Proposition}

\begin{document}

\title[Gradient Ricci soliton]
{Some properties of non-compact complete Riemannian manifolds}

\author{Li MA }

\address{Department of mathematical sciences \\
Tsinghua university \\
Beijing 100084 \\
China}

\email{lma@math.tsinghua.edu.cn} \dedicatory{}
\date{Nov. 20th, 2004}

\thanks{$^*$ This work is supported in part by
the Key 973 project of Ministry of Science and Technology of
China}

\begin{abstract}
In this paper, we study the volume growth property of a
non-compact complete Riemannian manifold $X$. We improve the
volume growth theorem of Calabi (1975) and Yau (1976), Cheeger,
Gromov and Taylor (1982). Then we use our new result to study
gradient Ricci solitons. We also show that on $X$, for any $q\in
(0,\infty)$, every non-negative $L^q$ subharmonic function is
constant under a natural decay condition on the Ricci curvature.
\end{abstract}

 \maketitle

\section{Introduction}
  The motivation for this paper comes from the interest in the
  understanding the Ricci solitons \cite{H95}. However, at this
  moment, almost all works in this direction are about Gradient
  Ricci Solitons. See \cite{H95},
\cite{B04},\cite{C96}, and \cite{M04}. Generally speaking, a
non-compact Ricci soliton may not be a gradient Ricci soliton. So
it may be interesting to consider problems related to Ricci
solitons.

  In this paper, we consider the volume growth properties
  of the non-compact complete Riemannian manifold $(X,g)$ under a natural Ricci
  curvature condition. We can improve the volume growth theorem
of Calabi \cite{Ca75} and Yau \cite{Y76}. Then we use our new
result to study gradient Ricci soliton. We also show that on $X$,
for any $q\in (0,\infty)$, every non-negative $L^q$ subharmonic
function is constant under a natural decay condition on the Ricci
curvature.

  Result about Gradient Solitons is stated in section four. In section two, we
  generalize the result of Calabi (\cite{Ca75}) and Yau (\cite{Y76})
  on infinite volume
  property for Riemannian manifolds with non-negative Ricci
  curvature. Calabi and Yau's result was generalized by Cheeger-Gromov-Taylor
  ( see Theorem 4.9
in \cite{CGT82} ) to
  Riemannian manifolds with lower bound like
  $$
Rc\geq -\frac{\nu_n}{r^2(x)},
  $$
for $r(x)>>1$ and some restricted dimensional constant $\nu_n$,
where $Rc$ is the Ricci tensor of the metric $g$, and $r(x)$ is
the distance function from some fixed point $x_0$. We can remove
this restriction to the dimensional constant $\nu_n$. Our result
is
\begin{Thm}
  Assume $(X,g)$ be a complete non-compact Riemannian manifold.
  Let $C(r)$ be a continuous function satisfying
$$
C(r)\geq -C{r^{-2}(x)}
$$
for $r$ large.
  Assume its Ricci curvature has the lower
bound
$$
Rc\geq C(r).
$$
Then for any $x\in X$ and $r>1$, there is a constant
$C(n,VolB_1(p))$ such that it holds
$$
VolB_r(x)\geq C(n,VolB_1(x))r.
$$
\end{Thm}

In section three. we give some remarks on Yau's gradient estimate
and vanishing properties for subharmonic functions.
\section{Proof of Theorem 1}

We let $D$ and $R(.,.).$ be the
 the Levi-Civita derivative and Riemmnain curvature of the metric $g$
 respectively. Let $x$ be a point
which is inside the cut locus of $p\in X $. Let $\gamma(t)$ be the
geodesic from $p$ to $x$. Choose a fixed point $x_0$ in the curve
$\gamma$. Let $r=r(x)$ be the distance function to the fixed point
$x_0\in X$ . For any $Y\in T_xX$ and $g(X,\frac{\partial}{\partial
t})=0$. Then we can get an Jacobi field $\hat{Y}$ by extending $X$
along $\gamma$ (see \cite{CE75} or \cite{Aub82}). Let $I_0^r(.,.)$
be the index form along $\gamma$. Then the hessian of $r$ at $x$
$$
H(r)(Y,Y)=\hat{Y}\hat{Y}r-D_{\hat{Y}}\hat{Y}r
$$
can be written as
$$
\int_0^r(|D_t\hat{Y}|^2-g(R(\hat{Y},D_t)D_t,\hat{Y}))dt
$$
which is the index form $I_0^r(\hat{Y},\hat{Y})$. We now extend
$Y$ along $\gamma$ and get a parallel vector field $E$. Then by
the minimizing property of the index form $I_0^r$ we have that
$$
I_0^r(\hat{Y},\hat{Y})\leq I_0^r(\frac{t}{r}E, \frac{t}{r}E),
$$
and the right side of the above inequality is
$$
\frac{1}{r}-\frac{1}{r^2}\int_0^r t^2g(R(E,D_t)D_t,E)dt
$$

To compute the Laplacian $\Delta r$ of $r$, we choose vector
fields
$$
\{ \frac{\partial}{\partial t}, E_1,.., E_{n-1} \}
$$
as an orthonormal basis of $T_{\gamma(t)}X$ and parallel along
$\gamma$. Then
 we have
 $$
\Delta r=\Sigma_{i=1}^{n-1}H(r)(E_i,E_i),
 $$
which is bounded above by
$$
\frac{n-1}{r}-\frac{1}{r^2}\int_0^r t^2Rc(D_t,D_t)dt.
$$
Using the assumption that $Rc\geq C(r)$ we get that
$$
\Delta r \leq \frac{n-1}{r}-\frac{1}{r^2}\int_0^r t^2C(t)dt.
$$
Since
$$
-\frac{1}{r^2}\int_0^r t^2C(t)dt\leq \frac{C}{r}.
$$
 we have
$$
\Delta r \leq \frac{n-1+C}{r}.
$$
Now it is standard to verify (see also page 7 in \cite{SY94} or
\cite{CM03} ) that in the distributional sense, it holds on $X$
$$
\Delta r \leq \frac{n-1+C}{r}.
$$
Then it holds in the distributional sense that
$$
\Delta r^2 =2r\Delta r+2\leq 2(n+1+C).
$$
That is, for any non-negative function $\phi\in C^{\infty}_0(X)$,
\begin{equation}
\int_X r^2\Delta \phi\leq 2(n+1+C)\int_X\phi.
\end{equation}
We now follow the argument of Schoen and Yau (see \cite{SY94}). By
approximation, we can let $\phi$ in (1) be a Lipschitz function
with compact support. Choose $\phi(x)=\xi(r(x))$, where $\xi(r)=1$
for $r\leq R-1$, $=0$ for $r\geq R+1$, and $\xi'(r)=-\frac{1}{2}$
for $R-1\leq r\leq R+1$. By direct computation we have
$$
\int_X r^2\Delta \phi=-2\int_{B_{R+1}(p)}\xi'r|Dr|^2
$$
Note that $|Dr|=1$, then we get
$$
\int_X r^2\Delta \phi= \int_{B_{R+1}(p)-B_{R-1}(p)}r
$$
From this we clearly have,
$$
\int_X r^2\Delta \phi\geq (R-1)Vol(B_{R+1}(p)-B_{R-1}(p)).
$$
Note that
$$
\int_X \phi \leq VolB_{R+1}(p)
$$
and
$$
B_1(x)\subset B_{R+1}(p)-B_{R-1}(p).
$$
Then by (1) we have
$$
(R-1)VolB_1(x)\leq 2(n+1+C)VolB_{R+1}(p).
$$
Since
$$
B_{R+1}(p)\subset B_{2(R+1)}(x),
$$
we obtain that
$$
(R-1)VolB_1(x)\leq 2(n+1+C)VolB_{2(R+1)}(x).
$$
This implies Theorem 1.

\section{Remarks on harmonic functions}

In his beautiful work \cite{Y75}, Yau proved that

\begin{Thm} Let $X$ be an $n$ ($\geq 2$) dimensional complete Riemannian
manifold with $Rc(X)\geq -(n-1)K$, where $K\geq 0$ is a constant.
Assume $u$ is a positive harmonic function on $X$. Let $B_R$ be a
geodesic ball in $X$. Then it holds on $B_{R/2}$ that
$$
|Dlog u|\leq C_n(\frac{1+R\sqrt{K}}{R})
$$
where $C_n$ is a constant depending only on $n$.
\end{Thm}

We now use this gradient estimate to study the positive solution
of the equation on the manifold $X$:
\begin{equation}
\Delta u=-c^2 u,
\end{equation}
where $c \geq 0$ is a constant.

\begin{Prop}
Let $X$ be an $n$ ($\geq 2$) dimensional complete Riemannian
manifold with $Rc(X)\geq -(n-1)K$, where $K\geq 0$ is a constant.
Assume $u$ is a positive solution of (2)  on $X$. Let $B_R$ be a
geodesic ball in $X$. Then it holds on $B_{R/2}$ that
$$
|Dlog u|\leq C_n(\frac{1+R\sqrt{K}}{R})
$$
where $C_n$ is a constant depending only on $n$.

\end{Prop}
\begin{proof} Let
$$
N=X\times R
$$
have the product metric, so that we still have $Rc(N)\geq -(n-1)K$
provided $Rc(X)\geq -(n-1)K$. We write $D_N$ as the covariant
derivative on $N$. Let
$$
w(x,t)=e^{ct}u(x).
$$
Then
$$
\Delta_N w=0.
$$
By Yau's estimate we have that
\begin{equation}
|Dlog u|=|Dlog w|\leq |D_Nlog w|\leq C_n(\frac{1+R\sqrt{K}}{R}).
\end{equation}
\end{proof}
The important matter for us is the constant $C_n$ independent of
the constant $c$. This is an important thing for us to study the
Ricci soliton.

We also need to study the positive solution of the following
equation on the manifold $X$:
\begin{equation}
\Delta u=c^2 u\geq 0,
\end{equation}
where $c \geq 0$ is a constant. Note that non-negative solutions
to (3) are non-negative subharmonic functions. So we now use our
Theorem 1 to study $L^q$ non-negative subharmonic functions on
complete Riemannian manifolds. As in the proof of Theorem 2.5 in
\cite{LS84}, we have

\begin{Prop} Suppose $X$ is a complete Riemannian manifold of
dimension $n$.  Let $C(r)$ be a continuous function satisfying
$$
C(r)\geq -C{r^{-2}(x)}
$$
for $r$ large.
  Assume that
there is a constant $C>0$ such that the Ricci curvature has the
bound
$$
Rc\geq C(r).
$$
Then for any $q\in (0, \infty)$, $X$ has no $L^q$ non-negative
subharmonic function except the constants.
\end{Prop}

\section{Ricci soliton with Ricci curvature quadratic decay}

In this section, we assume that the non-compact complete
Reimannian manifold $(X,g)$ is a gradient Ricci soliton, that is,
there is a smooth function $f$ such that
$$
Rc(g)=D^2f,
$$
on $X$. The classification for Ricci solitons is important to the
research for Ricci flow, see \cite{H95}, \cite{B04},\cite{C96},
\cite{Ma04}  and \cite{M04}. Generally speaking, a non-compact
Ricci soliton may not be a gradient Ricci soliton.

For a gradient Ricci soliton, as showed by R.Hamilton \cite{H95},
we have a constant $M$ such that
$$
|Df|^2+s=M
$$
where $s$ is the scalar curvature of $G$. We assume that $f$ is
not a constant, so $X$ is not Ricci flat.

We now come to the question: {\em Whether the constant $M$ is
bounded under a nice curvature condition?}

From the definition of the gradient Ricci soliton, it is clear
that
$$
s=\Delta f.
$$
Set
$$
u=e^f,
$$
which is a positive function on $X$. Then we have \cite{M04} that
$$
\Delta u=Mu,
$$
If $M$ is negative, then we can write it as $M=-c^2$, and then we
can use the gradient estimate in Proposition 3. Note that
$$
Dlog u=Df.$$ Then in this case, we have the bound
$$
|Df|\leq C_n(\frac{1+R\sqrt{K}}{R}),
$$
for every $R>0$, provided $Rc(X)\geq -(n-1)K$ on the ball $B_R$
and
$$ -n(n-1)K\leq M\leq 0.
$$

If $M=c^2$ for some constant $c\geq 0$ , we can use our
Proposition 4.

In conclusion we have
\begin{Thm}
Assume $(X,g)$ is a complete non-compact Gradient Ricci soliton
such that
$$
Rc=D^2f.
$$
Suppose $X$ is not Ricci flat.
  Let $C(r)$ be a continuous function satisfying
$$
C(r)\geq -C{r^{-2}(x)}
$$
for $r$ large.
  Assume the Ricci curvature has the
bound
$$
Rc\geq C(r).
$$
Then either (1) we have
$$
M=0,
$$
and $u:=e^f$ is a positive harmonic function on $X$;  or (2) $M>0$
and  for any $q\in (0, \infty)$, we have
$$
\int_X u^{q} =\infty.
$$
\end{Thm}

\begin{proof} In fact we have two cases when (i)
 $M\leq 0$ or (ii) $M> 0$. In case (i), we clearly have that
 $$
s<|Df|^2+s=M\leq 0.
 $$
The lower bound for $R$ follows easily from the assumption that
$$
s(x)\geq nC(r(x))
$$
for $r(x)$ large. By this we have $s(x)\to 0$ as $r(x)\to\infty$,
and $M=0$. In case (ii), $u$ is a positive subharmonic function,
amd we use Proposition 4 to get
$$
\int_X u^{q} =\infty.
$$
\end{proof}

\end{document}